\documentclass[12pt]{article}
\usepackage{amsthm}
\usepackage{amsmath}
\usepackage{amssymb}
\usepackage[utf8]{inputenc}
\usepackage{setspace}
\usepackage{relsize}
\newtheorem{theorem}{Theorem}[section]

\newtheorem{lemma}[theorem]{Lemma}

\newcommand{\davwebbbinom}{\genfrac{\langle}{\rangle}{0pt}{}}

\title{Generalized Lucas Theorem}
\author{Jordan Hirsh }
\date{January 2025}

\begin{document}

\maketitle

\begin{abstract}
    Let $p$ be a prime.  Let $A$ and $B$, $A \ge B \ge 0$, be integers with base $p$ expansions $A = \alpha_i\alpha_{i-1}\hdots \alpha_0$ and $B = \beta_i\beta_{i-1}\hdots \beta_0$.  Lucas proved that $$\binom{A}{B} \equiv \prod_{j=0}^{j=i}\binom{\alpha_j}{\beta_j} \text{ mod } p.$$ Similarly as proved by Kummer \cite{Kummer}, the $p$-adic valuation $v_p\binom{A}{B}$ is the number of borrows when computing $A-B$ in base $p$, or the number of carries in $(A-B)+B$ in base $p$. Davis and Webb \cite{DW} discovered a generalization of Lucas's Theorem for prime powers.  We prove a similar generalization in a different form using the concept of pseudo-digits.
\end{abstract}

\section{Introduction}

First we need to introduce the concept of pseudo-digits.  Let's start with an example $$\binom{432321433012}{323411244003}$$ in base $5$.  Starting at the right end, we have $2 < 3$ so we continue to the next digits, where we get $12 > 03$.  Thus $12,03$ is our first pair of pseudo-digits of $A,B$.  Now we go to the next digit and we have $0 \le 0$ so that is our second pseudo-digit.  We then have $3 < 4$, $33 < 44$, and finally $433 > 244$ which is our third pseudo-digit.  We keep appending digits to digits until we have a pseudo-digit $a_i$ greater than or equal to another pseudo-digit $b_i$. 
 Continuing in this way we get $$\binom{(4)(323)(2)(1)(433)(0)(12)}{(3)(234)(1)(1)(244)(0)(03)}.$$

 In more formal algorithmic terms, the way we determine the pseudo-digits $a_i$ is the following: We regard congruence classes mod $p^c$ as integers in $\{0,\hdots, p^c-1\}$ so we can order them.  Given $A$ and $B$ with $A \ge B$, we first determine the pseudo-digits $a_0$ and $b_0$.  We start at $c=1$ and see if $A$ mod $p^c \ge B$ mod $p^c$.  If not, we increase $c$ by $1$ and continue until we have $A$ mod $p^c \ge B$ mod $p^c$.  When that is true we say $a_0= A$ mod $p^c$ and likewise for $b_0$.  We then take $A' = \frac{(A-a_0)}{p^c}$ and $B' = \frac{(B-b_0)}{p^c}$ and we start the process again for $a_{1}$ and $b_{1}$ using $A'$ and $B'$.  We continue this process until $A'$ and $B'$ are both $0$, then we stop.  

 Note that the pseudo-digits of $A$ and $B$ depend on the pair $A,B$, not on $A$ and $B$ individually.

From now on in this paper Roman letters will refer to pseudo-digits whereas Greek letters will refer to singular digits.

If $a_i$ and $b_i$ are the $i$th pseudo-digits of $A$ and $B$ then $v_p\binom{a_i}{b_i}$ is one less than the number of base $p$ digits of $a_i$.  In other words, if $v_p \binom{a_j}{b_j} = k >0$ then $a_j$'s base $p$ representation is $\alpha_k\hdots\alpha_0$ and $b_j$'s is $\beta_k\hdots\beta_0$ where $\alpha_0 < \beta_0$, $\alpha_{\ell} \le \beta_{\ell}$ and $\alpha_k > \beta_k$.

Given a binomial coefficient $\binom{A}{B}$, we can write it as $p^mc$ where $p \nmid c$.  The following theorem gives what $c$ is mod $p^{n+m}$.  It is the main theorem of our paper.

\begin{theorem}
    Let $A \ge B \ge 0$ be integers and let $v_p\binom{A}{B} = m$.  Let $A = a_d\hdots a_0$ and $B = b_d\hdots b_0$ be the pseudo-digit expansions of $A$ and $B$.  Given $n \ge 1$,
    \begin{equation}\label{1}
        \binom{A}{B} = \binom{a_d\hdots a_0}{b_d\hdots b_0} \equiv \binom{a_d \hdots a_{d-n+1}}{b_d \hdots b_{d-n+1}} \prod_{i=0}^{d-n} \frac{\mathlarger{\binom{a_{i+n-1} \hdots a_i}{b_{i+n-1} \hdots b_i}}}{\mathlarger{\binom{a_{i+n-1} \hdots a_{i+1}}{b_{i+n-1} \hdots b_{i+1}}}} \mod  p^{n+m}.
    \end{equation} 
\end{theorem}

The numerator is a product of binomial coefficients formed from all blocks of $n$ consecutive pseudo-digits, whereas, except for omitting the initial factor, the denominator is the product of binomial coefficients formed from all blocks of $n-1$ consecutive pseudo-digits.  Note that the product is arranged so that the quotients are $p$ integral.  

Here is an example:  Take $$\binom{1221121202}{1011012021}$$ in base $3$.  Converting this to its pseudo-digit representation, we get $$\binom{(1)(2)(2)(1)(1)(21)(20)(2)}{(1)(0)(1)(1)(0)(12)(02)(1)}$$  We see there are 10 digits and 8 pseudo-digits so that means $v_3(\binom{1221121202}{1011012021}) = 10-8 = 2 = m$.

Suppose we want to find out what the above is mod $243$, which is $3^{2+3}$, so our $n$ is $3$.  Using the theorem, we get 
$$\binom{122}{101}\frac{\binom{221}{011}}{\binom{22}{01}}\frac{\binom{211}{110}}{\binom{21}{11}}\frac{\binom{1121}{1012}}{\binom{11}{10}}\frac{\binom{12120}{01202}}{\binom{121}{012}}\frac{\binom{21202}{12021}}{\binom{2120}{1202}} \text{ mod } 243.$$

To evaluate this expression, we use the following strategy: We write each factor as a power of $3$ times an integer mod $27$, where the latter integer is prime to $3$.  Since we know the total power of $3$ in the expression, we can write down the power of $3$ and multiply it by the product of mod $27$ parts.  
For example, consider the binomial coefficient $\binom{12120}{01202}$.  We know it's divisible by $9$ since there are two borrows.  So we evaluate it mod $243$ and we get $45$, which is $3^2 5$.  Now look at the denominator $\binom{121}{012}$.  This has one borrow, so we look at it mod $81$.  We obtain $75= 3^1 25$.  The quotient of these two is $\frac{3^2 5}{3^1 25} = \frac{3^1 5}{25} \equiv 3^1 11$ mod $81$.  The ratio contributes $11$ to the mod $27$ part of the final answer.  

Following through on all the factors we get
$$8\frac{14}{8}\frac{23}{8}\frac{30}{4}\frac{45}{75}\frac{90}{207} = $$
$$8\frac{14}{8}\frac{23}{8}\frac{3^1 10}{4}\frac{3^2 5}{3^1 25}\frac{3^2 10}{3^2 23} = $$
$$3^2 \: 8\,\frac{14}{8}\frac{23}{8}\frac{10}{4}\frac{5}{25}\frac{10}{23}$$

Simplifying this gets us $18$ mod $243$.  We see this is also correct by confirming with Pari-gp.

\section{Proof of Main Theorem}

Since Kummer showed that the $p$-adic valuation is the number of carries and pseudo-digits are defined based on when the carry operation ends, we have the following lemma:
\begin{lemma}\label{theremark}
    Consider $n$ consecutive pseudo-digits $a_{i+n-1} \hdots a_{i}$ and $b_{i+n-1} \hdots b_i$.  Then $$v_p\binom{a_{i+n-1} \hdots a_i}{b_{i+n-1} \hdots b_i} = -n+\sum_{j=i}^{i+n-1}\text{length }a_j  =\sum_{k=i}^{i+n-1}v_p\binom{a_k}{b_k}.$$  Further the $p$-adic valuation of $\frac{\binom{a_{i+n-1} \hdots a_i}{b_{i+n-1} \hdots b_i}}{\binom{a_{i+n-1} \hdots a_{i+1}}{b_{i+n-1} \hdots b_{i+1}}}$ is $$v_p\binom{a_i}{b_i}.$$
\end{lemma} 

We now start the proof of theorem.  The proof is by induction on $A+B$.

{\bf Case 0:} Base case

If we have an $A,B$ with less than or equal to $n$ pseudo-digits then it is trivially true that $\binom{A}{B} = \binom{a_{n-1}\hdots a_0}{b_{n-1}\hdots b_0}$, where the leading $a_i$'s are $0$ when the number of pseudo-digits is strictly less than $n$.  It is obvious that $\binom{a_{n-1}\hdots a_0}{b_{n-1}\hdots b_0} \equiv \binom{a_{n-1}\hdots a_0}{b_{n-1}\hdots b_0}$ mod $p^{n+m}$ where $m$ is the number of borrows. 

\medskip 

{\bf Case 1:}  $a_0 = 0$ and $b_0 = 0$.

Since $A = a_d\hdots a_1 0, B = b_d\hdots b_1 0$, we write $A=pA', B= pB'$ where $A' = a_d\hdots a_1, B' = b_d \hdots b_1$.  The first blocks of $n$ pseudo-digits of $A$ and $B$ can be written as $pA'' = a_{n-1}\hdots a_1 0$ and $pB''= b_{n-1}\hdots b_1 0$.   Since at minimum $pB''$ has $n$ pseudo-digits (we include the leading zeroes if there are any) and therefore at least $n$ digits, $B'-B'' \equiv 0$ mod $p^{n-1}$.  Let $k$ range from $B''+1$ to $B'$ and $j$ range over $1$ to $p-1$.
 Then $pk-j$ ranges over all relatively prime congruences classes mod $p^{n}$.  The same is true for $p(k+A'-B')-j$.  Thus the ratio of the respective products is $1$ mod $p^{n}$.  So we get 
\begin{equation}\label{2}
    \prod_{j=1}^{p-1}\prod_{k=B''+1}^{B'}\frac{p(k+A'-B')-j}{pk-j} \equiv 1 \text{ mod } p^{n}
\end{equation}
We then multiply both sides of $(\ref{2})$ by $$\prod_{j=1}^{p-1}\prod_{k=1}^{B''}\frac{p(k+A'-B')-j}{pk-j}$$ and get
$$\prod_{j=1}^{p-1}\prod_{k=1}^{B'}\frac{p(k+A'-B')-j}{pk-j} \equiv \prod_{j=1}^{p-1}\prod_{k=1}^{B''}\frac{p(k+A'-B')-j}{pk-j} \text{ mod } p^{n}$$

Since $pA'' \equiv pA'$ mod $p^n$ and $pB'' \equiv pB'$ mod $p^n$, we have $p(A''-B'') \equiv p(A'-B')$ mod $p^n$.  Thus the right side can be rewritten and we obtain

\begin{equation}\label{3}
    \prod_{j=1}^{p-1}\prod_{k=1}^{B'}\frac{p(k+A'-B')-j}{pk-j} \equiv \prod_{j=1}^{p-1}\prod_{k=1}^{B''}\frac{p(k+A''-B'')-j}{pk-j} \text{ mod } p^{n}
\end{equation}

As in Davis and Webb \cite{DW}, we rearrange factors as follows
\begin{gather*}
    \binom{A}{B} = \binom{pA'}{pB'} \\
    = \frac{\prod_{i=0}^{pB'-1}(pA'-i)}{\prod_{i=0}^{pB'-1}(pB'-i)} \\
    =\frac{\prod_{r=0}^{B'-1}p(A'-r)}{\prod_{r=0}^{B'-1}p(B'-r)}\prod_{j=1}^{p-1}\prod_{k=1}^{B'}\frac{p(k+A'-B')-j}{pk-j}\\
    = \frac{\prod_{r=0}^{B'-1}(A'-r)}{\prod_{r=0}^{B'-1}(B'-r)}\prod_{j=1}^{p-1}\prod_{k=1}^{B'}\frac{p(k+A'-B')-j}{pk-j} \\
    = \binom{A'}{B'}\prod_{j=1}^{p-1}\prod_{k=1}^{B'}\frac{p(k+A'-B')-j}{pk-j}
\end{gather*}

This allows us to rewrite both the left side and right side of  $(\ref{3})$ and obtain
$$\frac{\binom{A}{B}}{\binom{A'}{B'}} \equiv \frac{\binom{pA''}{pB''}}{\binom{A''}{B''}} \text{ mod } p^n$$
Multiplying both sides by $\binom{A'}{B'}$, we get 
\begin{equation}\label{4}
    \binom{A}{B} \equiv \binom{A'}{B'}\frac{\binom{pA''}{pB''}}{\binom{A''}{B''}} \text{ 
mod } p^{n+m}
\end{equation} since $v_{p}\binom{A'}{B'} = v_p\binom{A}{B} = m$.

We have
\begin{equation}\label{group1}
    \frac{\binom{pA''}{pB''}}{\binom{A''}{B''}} = \frac{\binom{a_{n-1}\hdots a_10}{b_{n-1}\hdots b_10}}{\binom{a_{n-1}\hdots a_1}{b_{n-1}\hdots b_1}}
\end{equation}
  and since $A'$ and $B'$ are shorter than $A$ and $B$, we invoke the induction hypothesis to obtain 

\begin{equation}\label{group2}
    \binom{A'}{B'} \equiv \binom{a_d \hdots a_{d-n+1}}{b_d \hdots b_{d-n+1}} \prod_{i=1}^{d-n} \frac{\mathlarger{\binom{a_{i+n-1} \hdots a_i}{b_{i+n-1} \hdots b_i}}}{\mathlarger{\binom{a_{i+n-1} \hdots a_{i+1}}{b_{i+n-1} \hdots b_{i+1}}}} \text{ mod } p^{n+m}
\end{equation}

Putting $(\ref{group1})$ and $(\ref{group2})$ together as in $(\ref{4})$ yields the theorem in this case:

\begin{equation}
    \binom{A}{B}\equiv\binom{a_d \hdots a_{d-n+1}}{b_d \hdots b_{d-n+1}} \prod_{i=0}^{d-n} \frac{\mathlarger{\binom{a_{i+n-1} \hdots a_i}{b_{i+n-1} \hdots b_i}}}{\mathlarger{\binom{a_{i+n-1} \hdots a_{i+1}}{b_{i+n-1} \hdots b_{i+1}}}} \text{ mod } p^{n+m}
\end{equation}

\medskip

{\bf Case 2:} $a_0 \not \equiv 0$ mod $p$ and $b_0 \not \equiv 0$ mod $p$.

If $A = B$ then the theorem is obvious, so we will assume $A > B$. 

Using the fact that $\binom{A}{B}B=\binom{A-1}{B-1}A$, we have
$$\binom{a_{n-1}\hdots a_0}{b_{n-1}\hdots b_0}(b_{n-1}\hdots b_0)=\binom{a_{n-1}\hdots a_0-1}{b_{n-1}\hdots b_0-1}(a_{n-1}\hdots a_0)$$
We divide by $\binom{a_{n-1}\hdots a_1}{b_{n-1}\hdots b_1}$ to get 
$$\frac{\mathlarger{\binom{a_{n-1}\hdots a_0}{b_{n-1}\hdots b_0}}}{\mathlarger{\binom{a_{n-1}\hdots a_1}{b_{n-1}\hdots b_1}}}(b_{n-1}\hdots b_0)=\frac{\mathlarger{\binom{a_{n-1}\hdots a_0-1}{b_{n-1}\hdots b_0-1}}}{\mathlarger{\binom{a_{n-1}\hdots a_1}{b_{n-1}\hdots b_1}}}(a_{n-1}\hdots a_0)$$
Since $(b_{n-1}\hdots b_0)$ and $(a_{n-1}\hdots a_0)$ are congruent to $B$ and $A$ mod $p^n$, we can invoke Lemma $\ref{theremark}$ to obtain $v_p\left(\frac{\binom{a_{n-1}\hdots a_0-1}{b_{n-1}\hdots b_0-1}}{\binom{a_{n-1}\hdots a_1}{b_{n-1}\hdots b_1}}\right) = v_p\binom{a_0-1}{b_0-1} \overset{\text{def}}{=} \ell$.  Similarly, $v_p\left(\frac{\binom{a_{n-1}\hdots a_0}{b_{n-1}\hdots b_0}}{\binom{a_{n-1}\hdots a_1}{b_{n-1}\hdots b_1}}\right) = v_p\binom{a_0}{b_0}$.   Since both $a_0$ and $b_0$ are non-zero, $v_p\binom{a_0}{b_0} = v_p\binom{a_0-1}{b_0-1}$.  Putting all this together, we have
\begin{equation}\label{stuff}
    \frac{\mathlarger{\binom{a_{n-1}\hdots a_0}{b_{n-1}\hdots b_0}}}{\mathlarger{\binom{a_{n-1}\hdots a_1}{b_{n-1}\hdots b_1}}}B\equiv\frac{\mathlarger{\binom{a_{n-1}\hdots a_0-1}{b_{n-1}\hdots b_0-1}}}{\mathlarger{\binom{a_{n-1}\hdots a_1}{b_{n-1}\hdots b_1}}}A \text{ mod } p^{n+\ell}
\end{equation}

Again invoking Lemma $\ref{theremark}$, we see

\begin{equation}\label{stuff2}
    \begin{gathered}
    v_p\left(\binom{a_d \hdots a_{d-n+1}}{b_d \hdots b_{d-n+1}} \prod_{i=1}^{d-n} \frac{\mathlarger{\binom{a_{i+n-1} \hdots a_i}{b_{i+n-1} \hdots b_i}}}{\mathlarger{\binom{a_{i+n-1} \hdots a_{i+1}}{b_{i+n-1} \hdots b_{i+1}}}}\right) \\
    =\sum_{i=1}^{d} v_p\binom{a_i}{b_i} = \sum_{i=0}^{d} v_p\binom{a_i}{b_i}-v_p\binom{a_0}{b_0} = m-\ell. 
\end{gathered}
\end{equation}

Now multiplying both sides of $(\ref{stuff})$ by the top expression in $(\ref{stuff2})$, we obtain
\begin{equation}\label{6}
    \begin{gathered}
        \binom{a_d \hdots a_{d-n+1}}{b_d \hdots b_{d-n+1}} \left(\prod_{i=1}^{d-n} \frac{\mathlarger{\binom{a_{i+n-1} \hdots a_i}{b_{i+n-1} \hdots b_i}}}{\mathlarger{\binom{a_{i+n-1} \hdots a_{i+1}}{b_{i+n-1} \hdots b_{i+1}}}}\right)\frac{\mathlarger{\binom{a_{n-1}\hdots a_0}{b_{n-1}\hdots b_0}}}{\mathlarger{\binom{a_{n-1}\hdots a_1}{b_{n-1}\hdots b_1}}}B \equiv \\
        \binom{a_d \hdots a_{d-n+1}}{b_d \hdots b_{d-n+1}} \left(\prod_{i=1}^{d-n} \frac{\mathlarger{\binom{a_{i+n-1} \hdots a_i}{b_{i+n-1} \hdots b_i}}}{\mathlarger{\binom{a_{i+n-1} \hdots a_{i+1}}{b_{i+n-1} \hdots b_{i+1}}}}\right)\frac{\mathlarger{\binom{a_{n-1}\hdots a_0-1}{b_{n-1}\hdots b_0-1}}}{\mathlarger{\binom{a_{n-1}\hdots a_1}{b_{n-1}\hdots b_1}}}A \text{ mod } p^{n+m}
    \end{gathered}
\end{equation}

By the induction hypothesis, the second half of $(\ref{6})$ is congruent to $\binom{A-1}{B-1}A$.  Since we have $\binom{A}{B}B = \binom{A-1}{B-1}A$, $(\ref{6})$ becomes 
$$\binom{a_d \hdots a_{d-n+1}}{b_d \hdots b_{d-n+1}} \left(\prod_{i=1}^{d-n} \frac{\mathlarger{\binom{a_{i+n-1} \hdots a_i}{b_{i+n-1} \hdots b_i}}}{\mathlarger{\binom{a_{i+n-1} \hdots a_{i+1}}{b_{i+n-1} \hdots b_{i+1}}}}\right)\frac{\mathlarger{\binom{a_{n-1}\hdots a_0}{b_{n-1}\hdots b_0}}}{\mathlarger{\binom{a_{n-1}\hdots a_1}{b_{n-1}\hdots b_1}}}B \equiv \binom{A}{B}B \text{ mod } p^{n+m}$$
Since $p$$\not|B$, this is equivalent to
$$\binom{A}{B} \equiv \binom{a_d \hdots a_{d-n+1}}{b_d \hdots b_{d-n+1}} \prod_{i=0}^{d-n} \frac{\mathlarger{\binom{a_{i+n-1} \hdots a_i}{b_{i+n-1} \hdots b_i}}}{\mathlarger{\binom{a_{i+n-1} \hdots a_{i+1}}{b_{i+n-1} \hdots b_{i+1}}}} \text{ mod } p^{n+m}$$

{\bf Case 3}: $a_0 \not \equiv 0$ mod $p$ and $b_0 = 0$.

Using the fact that $\binom{A}{B}(A-B) = \binom{A-1}{B}A$, we have

$$\binom{a_{n-1}\hdots a_0}{b_{n-1}\hdots b_0}(a_n\hdots a_0-b_n\hdots b_0)=\binom{a_{n-1}\hdots a_0-1}{b_{n_1}\hdots b_0}(a_{n-1}\hdots a_0)$$ Dividing  gets us

$$\frac{\mathlarger{\binom{a_{n-1}\hdots a_0}{b_{n-1}\hdots b_0}}}{\mathlarger{\binom{a_{n-1}\hdots a_1}{b_{n-1}\hdots b_1}}}(a_{n-1}\hdots a_0-b_{n-1}\hdots b_0)=\frac{\mathlarger{\binom{a_{n-1}\hdots a_0-1}{b_{n-1}\hdots b_0}}}{\mathlarger{\binom{a_{n-1}\hdots a_1}{b_{n-1}\hdots b_1}}}(a_{n-1}\hdots a_0)$$ and since $(a_{n-1}\hdots a_0)$ and $(b_{n-1}\hdots b_0)$ are congruent to $A$ and $B$ mod $p^n$, respectively, the above simplifies to

\begin{equation}\label{thequation}
    \frac{\mathlarger{\binom{a_{n-1}\hdots a_0}{b_{n-1}\hdots b_0}}}{\mathlarger{\binom{a_{n-1}\hdots a_1}{b_{n-1}\hdots b_1}}}(A-B)\equiv \frac{\mathlarger{\binom{a_{n-1}\hdots a_0-1}{b_{n-1}\hdots b_0}}}{\mathlarger{\binom{a_{n-1}\hdots a_1}{b_{n-1}\hdots b_1}}}A \text{ mod } p^{n}
\end{equation}
As in Case 2, we multiply the above expression by
\begin{equation}\label{thequation2}
    \binom{a_d \hdots a_{d-n+1}}{b_d \hdots b_{d-n+1}} \prod_{i=1}^{d-n} \frac{\mathlarger{\binom{a_{i+n-1} \hdots a_i}{b_{i+n-1} \hdots b_i}}}{\mathlarger{\binom{a_{i+n-1} \hdots a_{i+1}}{b_{i+n-1} \hdots b_{i+1}}}}
\end{equation}
Again we invoke Lemma $\ref{theremark}$ to find the $p$-adic valuation of $(\ref{thequation2})$ is

$$\sum_{i=1}^{d} v_p\binom{a_i}{b_i} = \sum_{i=0}^{d} v_p\binom{a_i}{b_i}-v_p\binom{a_0}{b_0}.$$

However $v_p\binom{a_0}{b_0} = 0$ so the above expression is equal to $v_p\binom{A}{B}=m$.  Thus when mulitiplying $(\ref{thequation})$ by $(\ref{thequation2})$ we get

\begin{equation}\label{7}
    \begin{gathered}
        \binom{a_d \hdots a_{d-n+1}}{b_d \hdots b_{d-n+1}} \left(\prod_{i=1}^{d-n} \frac{\mathlarger{\binom{a_{i+n-1} \hdots a_i}{b_{i+n-1} \hdots b_i}}}{\mathlarger{\binom{a_{i+n-1} \hdots a_{i+1}}{b_{i+n-1} \hdots b_{i+1}}}}\right)\frac{\mathlarger{\binom{a_{n-1}\hdots a_0}{b_{n-1}\hdots b_0}}}{\mathlarger{\binom{a_{n-1}\hdots a_1}{b_{n-1}\hdots b_1}}}(A-B)\equiv\\
        \binom{a_d \hdots a_{d-n+1}}{b_d \hdots b_{d-n+1}} \left(\prod_{i=1}^{d-n} \frac{\mathlarger{\binom{a_{i+n-1} \hdots a_i}{b_{i+n-1} \hdots b_i}}}{\mathlarger{\binom{a_{i+n-1} \hdots a_{i+1}}{b_{i+n-1} \hdots b_{i+1}}}}\right)\frac{\mathlarger{\binom{a_{n-1}\hdots a_0-1}{b_{n-1}\hdots b_0}}}{\mathlarger{\binom{a_{n-1}\hdots a_1}{b_{n-1}\hdots b_1}}}A \text{ mod } p^{n+m}
    \end{gathered}
\end{equation}

By the induction hypothesis, the second expression is congruent to $\binom{A-1}{B}A$ mod $p^{n+m}$. Since $\binom{A}{B}(A-B)=\binom{A-1}{B}A$, the rest follows as in the end of Case 2.

{\bf Case 4}: $a_0 = 0$ and $b_0 \not \equiv 0$ mod $p$.

We start with the identity $\binom{A}{B}B = \binom{A}{B-1}(A-B+1)$.  There are two subcases, Case $4a$: when $B \not \equiv 1$ mod $p$, and Case $4b$: when $B \equiv 1$ mod $p$.

{\bf Case 4a}: $B \not \equiv 1$ mod $p$:

Both $B$ and $A-B+1$ are $\not \equiv 0$ mod $p$.  We have
$$\binom{a_{n-1}\hdots a_0}{b_{n-1}\hdots b_0}(b_{n-1}\hdots b_0) = \binom{a_{n-1}\hdots a_0}{b_{n-1}\hdots b_0-1}(a_{n-1}\hdots a_0-b_{n-1}\hdots b_0+1)$$

Dividing by $\binom{a_{n-1}\hdots a_1}{b_{n-1}\hdots b_1}$ yields $$\frac{\mathlarger{\binom{a_{n-1}\hdots a_0}{b_{n-1}\hdots b_0}}}{\mathlarger{\binom{a_{n-1}\hdots a_1}{b_{n-1}\hdots b_1}}}(b_{n-1}\hdots b_0)=\frac{\mathlarger{\binom{a_{n-1}\hdots a_0}{b_{n-1}\hdots b_0-1}}}{\mathlarger{\binom{a_{n-1}\hdots a_1}{b_{n-1}\hdots b_1}}}(a_{n-1}\hdots a_0-b_{n-1}\hdots b_0+1)$$ and since $a_{n-1}\hdots a_0$ and $b_{n-1}\hdots b_0$ are congruent to $A$ and $B$ mod $p^n$, respectively, by invoking Lemma $\ref{theremark}$, we have

\begin{equation}\label{filler}
    \frac{\mathlarger{\binom{a_{n-1}\hdots a_0}{b_{n-1}\hdots b_0}}}{\mathlarger{\binom{a_{n-1}\hdots a_1}{b_{n-1}\hdots b_1}}}B\equiv\frac{\mathlarger{\binom{a_{n-1}\hdots a_0}{b_{n-1}\hdots b_0-1}}}{\mathlarger{\binom{a_{n-1}\hdots a_1}{b_{n-1}\hdots b_1}}}(A-B+1) \text{ mod } p^{n+\ell}
\end{equation}
where $\ell$ as in Case $2$ is $v_p\binom{a_0}{b_0} = v_p\binom{a_0}{b_0-1}$.

As usual, we then multiply by

$$\binom{a_d \hdots a_{d-n+1}}{b_d \hdots b_{d-n+1}} \prod_{i=1}^{d-n} \frac{\mathlarger{\binom{a_{i+n-1} \hdots a_i}{b_{i+n-1} \hdots b_i}}}{\mathlarger{\binom{a_{i+n-1} \hdots a_{i+1}}{b_{i+n-1} \hdots b_{i+1}}}}$$ which as in Case 1 has $p$-adic valuation $m-\ell$.  Therefore,

\begin{equation}\label{8}
    \begin{gathered}
        \binom{a_d \hdots a_{d-n+1}}{b_d \hdots b_{d-n+1}} \left(\prod_{i=1}^{d-n} \frac{\mathlarger{\binom{a_{i+n-1} \hdots a_i}{b_{i+n-1} \hdots b_i}}}{\mathlarger{\binom{a_{i+n-1} \hdots a_{i+1}}{b_{i+n-1} \hdots b_{i+1}}}}\right)\frac{\mathlarger{\binom{a_{n-1}\hdots a_0}{b_{n-1}\hdots b_0}}}{\mathlarger{\binom{a_{n-1}\hdots a_1}{b_{n-1}\hdots b_1}}}B\equiv\\
        \binom{a_d \hdots a_{d-n+1}}{b_d \hdots b_{d-n+1}} \left(\prod_{i=1}^{d-n} \frac{\mathlarger{\binom{a_{i+n-1} \hdots a_i}{b_{i+n-1} \hdots b_i}}}{\mathlarger{\binom{a_{i+n-1} \hdots a_{i+1}}{b_{i+n-1} \hdots b_{i+1}}}}\right)\frac{\mathlarger{\binom{a_{n-1}\hdots a_0}{b_{n-1}\hdots b_0-1}}}{\mathlarger{\binom{a_{n-1}\hdots a_1}{b_{n-1}\hdots b_1}}}(A-B+1) \\ \text{ mod } p^{n+m}
    \end{gathered}
\end{equation}

By invoking the induction hypothesis on the bottom part of $(\ref{8})$, we see that it is congruent to $\binom{A}{B-1}(A-B+1)$ mod $p^{n+m}$.  Since $\binom{A}{B}B = \binom{A}{B-1}(A-B+1)$, the rest follows as in the end of Case 2. 

{\bf Case 4b}: $B \equiv 1$ mod $p$. 

We again have $\binom{A}{B}B = \binom{A}{B-1}(A-B+1)$.  In this case $\binom{A}{B-1}$ is of the form $\binom{a_d\hdots a_0}{b_d \hdots b_0-1}$.  We have the base $p$ expansions $a_0 = \alpha_{\ell-1}\hdots\alpha_10$ and $b_0-1 = \beta_{\ell-1}\hdots\beta_10$.  Since $a_0,b_0$ is a pair of pseudo-digits, $\alpha_{\ell-1} > \beta_{\ell-1}$ and $\alpha_k \le \beta_k$ for $1 \le k < \ell-1$.  Since the $p$-adic expansions both end in $0$, these $0$'s give us a pair of pseudo-digits for $a_0,b_0-1$.  Now we have new pseudo-digits ending in $\alpha_1$ and $\beta_1$.  If $\alpha_1 = \beta_1$ then this is also a pair of pseudo-digits.  We continue breaking off one-digit pairs of pseudo-digits until either $\alpha_f < \beta_f$ for some $f < \ell-1$ or we get to $\alpha_{\ell-1} > \beta_{\ell-1}$, in which case we let $f=\ell-1$.  The pseudo-digit expansions of $a_{n-1}\hdots a_0$ and $b_{n-1}\hdots b_0-1$ are $ a_{n-1}\hdots a_0'\alpha_{f-1}\hdots\alpha_10$ and $b_{n-1}\hdots b_0'\beta_{f-1}\hdots\beta_10$ where $a_0' = \alpha_{\ell-1}\hdots\alpha_f$ and $b_0' = \beta_{\ell-1}\hdots\beta_f$.  As before, we have
\begin{equation}\label{something}
    \frac{\mathlarger{\binom{a_{n-1}\hdots a_0}{b_{n-1}\hdots b_0}}}{\mathlarger{\binom{a_{n-1}\hdots a_1}{b_{n-1}\hdots b_1}}}(b_{n-1}\hdots b_0)=\frac{\mathlarger{\binom{a_{n-1}\hdots a_0}{b_{n-1}\hdots b_0-1}}}{\mathlarger{\binom{a_{n-1}\hdots a_1}{b_{n-1}\hdots b_1}}}(a_{n-1}\hdots a_0-b_{n-1}\hdots b_0+1)
\end{equation}

Since $v_p(\frac{\binom{a_{n-1}\hdots a_0}{b_{n-1}\hdots b_0}}{\binom{a_{n-1}\hdots a_1}{b_{n-1}\hdots b_1}}) = v_p\binom{a_0}{b_0} = \ell$ and $B \equiv b_{n-1}\hdots b_0 \mod p^n$, we have $$\frac{\mathlarger{\binom{a_{n-1}\hdots a_0}{b_{n-1}\hdots b_0}}}{\mathlarger{\binom{a_{n-1}\hdots a_1}{b_{n-1}\hdots b_1}}}(b_{n-1}\hdots b_0) \equiv \frac{\mathlarger{\binom{a_{n-1}\hdots a_0}{b_{n-1}\hdots b_0}}}{\mathlarger{\binom{a_{n-1}\hdots a_1}{b_{n-1}\hdots b_1}}}B \text{ mod } p^{n+\ell}.$$
To handle the right side of $(\ref{something})$, first note that 
$\frac{\binom{a_{n-1}\hdots a_0}{b_{n-1}\hdots b_0-1}}{\binom{a_{n-1}\hdots a_1}{b_{n-1}\hdots b_1}}=\frac{\binom{a_{n-1}\hdots a_0'\alpha_{f-1}\hdots \alpha_1 0}{b_{n-1}\hdots b_0' \beta_{f-1}\hdots \beta_1 0}}{\binom{a_{n-1}\hdots a_1}{b_{n-1}\hdots b_1}}$.  Moreover, by Lemma $\ref{theremark}$,
\begin{equation}\label{howmanymoreofthesewillihave}
    \begin{split}
    &v_p\left(\frac{\mathlarger{\binom{a_{n-1}\hdots a_0'\alpha_{f-1}\hdots \alpha_1 0}{b_{n-1}\hdots b_0' \beta_{f-1}\hdots \beta_1 0}}}{\mathlarger{\binom{a_{n-1}\hdots a_1}{b_{n-1}\hdots b_1}}}\right) \\ 
    &= \sum^{n-1}_{i=1} v_p\binom{a_i}{b_i}+v_p\binom{a_0'}{b_0'}+\sum^{f-1}_{j=0}v_p\binom{\alpha_j}{\beta_j} - \sum^{n-1}_{i=1} v_p\binom{a_i}{b_i} \\
    &= v_p\binom{a_0'}{b_0'}+\sum^{f-1}_{j=0}v_p\binom{\alpha_j}{\beta_j} = \ell-f+\sum^{f-1}_{j=0} 0 = \ell-f.
    \end{split}
\end{equation}

Because $a_{n-1}\hdots a_0'\alpha_{f-1}\hdots\alpha_10$ and $ b_{n-1}\hdots b_0'\beta_{f-1}\hdots\beta_10$ agree with $A$ and $ B-1$ on at least the first $n+f$ digits we have 
$$a_{n-1}\hdots a_0-b_{n-1}\hdots b_0+1 \equiv A-B+1 \text{ mod } p^{n+f}$$ therefore

\begin{equation}\label{nearingtheend}
\begin{gathered}
    \frac{\mathlarger{\binom{a_{n-1}\hdots a_0}{b_{n-1}\hdots b_0-1}}}{\mathlarger{\binom{a_{n-1}\hdots a_1}{b_{n-1}\hdots b_1}}}(a_{n-1}\hdots a_0-b_{n-1}\hdots b_0+1) \equiv \\ \frac{\mathlarger{\binom{a_{n-1}\hdots a_0}{b_{n-1}\hdots b_0-1}}}{\mathlarger{\binom{a_{n-1}\hdots a_1}{b_{n-1}\hdots b_1}}}(A-B+1) \text{ mod } p^{n+\ell}
\end{gathered}
\end{equation}

So applying all this to $(\ref{something})$ gives us 

\begin{equation}\label{9}
       \frac{\mathlarger{\binom{a_{n-1}\hdots a_0}{b_{n-1}\hdots b_0}}}{\mathlarger{\binom{a_{n-1}\hdots a_1}{b_{n-1}\hdots b_1}}}B\equiv \frac{\mathlarger{\binom{a_{n-1}\hdots a_0'\alpha_{f-1}\hdots\alpha_10}{b_{n-1}\hdots b_0'\beta_{f-1}\hdots\beta_10}}}{\mathlarger{\binom{a_{n-1}\hdots a_1}{b_{n-1}\hdots b_1}}}(A-B+1)\text{ mod } p^{n+\ell}
\end{equation}

The $p$-adic valuation of $\binom{a_{n-1}\hdots a_0'\alpha_{f-1}\hdots\alpha_10}{b_{n-1}\hdots b_0'\beta_{f-1}\hdots\beta_10}$ is the number of borrows, which is the number of borrows in $a_{n-1}\hdots a_0' - b_{n-1} \hdots b_0'$ since $\alpha_j = \beta_j$ for $j < f$.  So the $p$-adic valuation is $s-n$ where $s$ is the number of $p$-adic digits of $a_{n-1}\hdots a_0'$.  Invoking the induction hypothesis, we get

$$\binom{a_{n-1}\hdots a_0'\alpha_{f-1}\hdots\alpha_10}{b_{n-1}\hdots b_0'\beta_{f-1}\hdots\beta_10} \equiv \binom{a_{n-1} \hdots a_{0}'}{b_{n-1} \hdots b_{0}'} \prod_{j=0}^{f-1} \frac{\mathlarger{\binom{\alpha_{j+n-1} \hdots \alpha_j}{\beta_{j+n-1} \hdots \beta_j}}}{\mathlarger{\binom{\alpha_{j+n-1} \hdots \alpha_{j+1}}{\beta_{j+n-1} \hdots \beta_{j+1}}}} \text{ mod } p^{n+(s-n)}$$
where we denote $\alpha_f = a_0', \beta_f = b_0'$ and $\alpha_{f+i}=a_i, \beta_{f+i}=b_i$.  Since $\alpha_{j} = \beta_j$ for $j < f$, we have $v_p\left(\frac{\binom{\alpha_{j+n-1} \hdots \alpha_j}{\beta_{j+n-1} \hdots \beta_j}}{\binom{\alpha_{j+n-1} \hdots \alpha_{j+1}}{\beta_{j+n-1} \hdots \beta_{j+1}}}\right)=0$ for each $j < f$.

We must divide by $\binom{a_{n-1}\hdots a_1}{b_{n-1}\hdots b_1}$ which has valuation $(s-n-\ell+f)$.  This is demonstrated by invoking $(\ref{howmanymoreofthesewillihave})$ and we get:
$$v_p\binom{a_{n-1}\hdots a_1}{b_{n-1}\hdots b_1} = \sum_{i=1}^{n-1}v_p\binom{a_i}{b_i} = \sum_{i=1}^{n-1}v_p\binom{a_i}{b_i}+v_p\binom{a_0'}{b_0'}-v_p\binom{a_0'}{b_0'}=$$ $$s-n-(\ell-f).$$
This then gives us
$$\frac{\mathlarger{\binom{a_{n-1}\hdots a_0'\alpha_{f-1}\hdots\alpha_10}{b_{n-1}\hdots b_0'\beta_{f-1}\hdots\beta_10}}}{\mathlarger{\binom{a_{n-1}\hdots a_1}{b_{n-1}\hdots b_1}}} \equiv \frac{\mathlarger{\binom{a_{n-1} \hdots a_{0}'}{b_{n-1} \hdots b_{0}'}}}{\mathlarger{\binom{a_{n-1}\hdots a_1}{b_{n-1}\hdots b_1}}} \prod_{j=0}^{f-1} \frac{\mathlarger{\binom{\alpha_{j+n-1} \hdots \alpha_j}{\beta_{j+n-1} \hdots \beta_j}}}{\mathlarger{\binom{\alpha_{j+n-1} \hdots \alpha_{j+1}}{\beta_{j+n-1} \hdots \beta_{j+1}}}} \text{ mod } p^{n+(s-n)-(s-n-\ell+f)}$$
The above expression simplifies to
$$\frac{\mathlarger{\binom{a_{n-1}\hdots a_0'\alpha_{f-1}\hdots\alpha_10}{b_{n-1}\hdots b_0'\beta_{f-1}\hdots\beta_10}}}{\mathlarger{\binom{a_{n-1}\hdots a_1}{b_{n-1}\hdots b_1}}} \equiv \frac{\mathlarger{\binom{a_{n-1} \hdots a_{0}'}{b_{n-1} \hdots b_{0}'}}}{\mathlarger{\binom{a_{n-1}\hdots a_1}{b_{n-1}\hdots b_1}}} \prod_{j=0}^{f-1} \frac{\mathlarger{\binom{\alpha_{j+n-1} \hdots \alpha_j}{\beta_{j+n-1} \hdots \beta_j}}}{\mathlarger{\binom{\alpha_{j+n-1} \hdots \alpha_{j+1}}{\beta_{j+n-1} \hdots \beta_{j+1}}}} \text{ mod } p^{n+\ell-f}$$
So taking this altogether $(\ref{9})$ is

\begin{equation}\label{thisisalot}
    \frac{\mathlarger{\binom{a_{n-1}\hdots a_0}{b_{n-1}\hdots b_0}}}{\mathlarger{\binom{a_{n-1}\hdots a_1}{b_{n-1}\hdots b_1}}}B\equiv \frac{\mathlarger{\binom{a_{n-1} \hdots a_{0}'}{b_{n-1} \hdots b_{0}'}}}{\mathlarger{\binom{a_{n-1}\hdots a_1}{b_{n-1}\hdots b_1}}} \prod_{j=0}^{f-1} \frac{\mathlarger{\binom{\alpha_{j+n-1} \hdots \alpha_j}{\beta_{j+n-1} \hdots \beta_j}}}{\mathlarger{\binom{\alpha_{j+n-1} \hdots \alpha_{j+1}}{\beta_{j+n-1} \hdots \beta_{j+1}}}} (A-B+1)\text{ mod } p^{n+\ell}
\end{equation}
We then multiply by 

$$\binom{a_d \hdots a_{d-n+1}}{b_d \hdots b_{d-n+1}} \prod_{i=1}^{d-n} \frac{\mathlarger{\binom{a_{i+n-1} \hdots a_i}{b_{i+n-1} \hdots b_i}}}{\mathlarger{\binom{a_{i+n-1} \hdots a_{i+1}}{b_{i+n-1} \hdots b_{i+1}}}}$$
which as in Case 1, has valuation $m-\ell$.  

This yields

\begin{equation}\label{10}
    \begin{gathered}
        \binom{a_d \hdots a_{d-n+1}}{b_d \hdots b_{d-n+1}} \prod_{i=1}^{d-n} \frac{\mathlarger{\binom{a_{i+n-1} \hdots a_i}{b_{i+n-1} \hdots b_i}}}{\mathlarger{\binom{a_{i+n-1} \hdots a_{i+1}}{b_{i+n-1} \hdots b_{i+1}}}}\frac{\mathlarger{\binom{a_{n-1}\hdots a_0}{b_{n-1}\hdots b_0}}}{\mathlarger{\binom{a_{n-1}\hdots a_1}{b_{n-1}\hdots b_1}}}B\equiv\\
        \binom{a_d \hdots a_{d-n+1}}{b_d \hdots b_{d-n+1}} \prod_{i=1}^{d-n} \frac{\mathlarger{\binom{a_{i+n-1} \hdots a_i}{b_{i+n-1} \hdots b_i}}}{\mathlarger{\binom{a_{i+n-1} \hdots a_{i+1}}{b_{i+n-1} \hdots b_{i+1}}}}\frac{\mathlarger{\binom{a_{n-1}\hdots a_0'}{b_{n-1}\hdots b_0'}}}{\mathlarger{\binom{a_{n-1}\hdots a_1}{b_{n-1}\hdots b_1}}}\prod_{j=0}^{f-1}\frac{\mathlarger{\binom{\alpha_{j+n-1} \hdots \alpha_j}{\beta_{j+n-1} \hdots \beta_j}}}{\mathlarger{\binom{\alpha_{j+n-1} \hdots \alpha_{j+1}}{\beta_{j+n-1} \hdots \beta_{j+1}}}}(A-B+1)\\ 
        \text{ mod } p^{n+m}
    \end{gathered}
\end{equation}

By the induction hypothesis the middle expression is equivalent to $\binom{A}{B-1}(A-B+1)$ and we have $\binom{A}{B}B = \binom{A}{B-1}(A-B+1)$ which means
$$\binom{A}{B} \equiv \binom{a_d \hdots a_{d-n+1}}{b_d \hdots b_{d-n+1}} \prod_{i=0}^{d-n} \frac{\binom{a_{i+n-1} \hdots a_i}{b_{i+n-1} \hdots b_i}}{\binom{a_{i+n-1} \hdots a_{i+1}}{b_{i+n-1} \hdots b_{i+1}}} \text{ mod } p^{n+m}.$$

\section{Comparison}
We now compare Davis and Webb's version to our version.  Let $A = \alpha_d\alpha_{d-1}\hdots\alpha_0$ and $B=\beta_d\beta_{d-1}\hdots\beta_0$ in base $p$.  

Define $\davwebbbinom{\alpha}{\beta} = p$ if $\alpha < \beta$ and $$\davwebbbinom{\alpha_{i+N-1}\hdots\alpha_i}{\beta_{i+N-1}\hdots\beta_i} = 
\begin{cases} 
\binom{\mathlarger{\alpha_{i+N-1}\hdots\alpha_i}}{\mathlarger{\beta_{i+N-1}\hdots\beta_i}} \text{ if } \alpha_{i+N-1}\hdots\alpha_i > \beta_{i+N-1}\hdots\beta_i \\ \\
p\davwebbbinom{\mathlarger{\alpha_{i+N-2}\hdots\alpha_i}}{\mathlarger{\beta_{i+N-2}\hdots\beta_i}} \text{ otherwise.}
\end{cases}$$

Davis and Webb prove $$\binom{A}{B} \equiv \davwebbbinom{\alpha_d\hdots\alpha_{d-N+1}}{\beta_d\hdots\beta_{d-N+1}}\prod_{i=0}^{d-N}\frac{\mathlarger{\davwebbbinom{\alpha_{i+N-1}\hdots\alpha_i}{\beta_{i+N-1}\hdots\beta_i}}}{\mathlarger{\davwebbbinom{\alpha_{i+N-1}\hdots\alpha_{i+1}}{\beta_{i+N-1}\hdots\beta_{i+1}}}} \text{ mod }  p^N.$$

We proved $$\binom{A}{B}\equiv\binom{a_d \hdots a_{d-n+1}}{b_d \hdots b_{d-n+1}} \prod_{i=0}^{d-n} \frac{\mathlarger{\binom{a_{i+n-1} \hdots a_i}{b_{i+n-1} \hdots b_i}}}{\mathlarger{\binom{a_{i+n-1} \hdots a_{i+1}}{b_{i+n-1} \hdots b_{i+1}}}} \mod  p^{n+m}$$ where $m$ is the number of borrows in the expression and $n+m$ is the exponent of the modulus we are looking at.

As an example take $\binom{21202112}{12021110}$ and look at it mod $243$.  Davis and Webb would put it in this form:
$$\davwebbbinom{21202}{12021}\frac{\mathlarger{\davwebbbinom{12021}{20211}}}{\mathlarger{\davwebbbinom{1202}{2021}}}\frac{\mathlarger{\davwebbbinom{20211}{02111}}}{\mathlarger{\davwebbbinom{2021}{0211}}}\frac{\mathlarger{\davwebbbinom{02112}{21110}}}{\mathlarger{\davwebbbinom{0211}{2111}}} \text{ mod } 243$$  which converts into $$\binom{21202}{12021}\frac{3\mathlarger{\binom{2021}{0211}}}{3\mathlarger{\binom{202}{021}}}\frac{\mathlarger{\binom{20211}{02111}}}{\mathlarger{\binom{2021}{0211}}}\frac{3\mathlarger{\binom{2112}{1110}}}{3\mathlarger{\binom{211}{111}}} \equiv \binom{21202}{12021}\frac{\mathlarger{\binom{2021}{0211}}}{\mathlarger{\binom{202}{021}}}\frac{\mathlarger{\binom{20211}{02111}}}{\mathlarger{\binom{2021}{0211}}}\frac{\mathlarger{\binom{2112}{1110}}}{\mathlarger{\binom{211}{111}}} \text{ mod } 243$$ which is $90\frac{39}{3}\frac{66}{39}\frac{152}{242} \equiv 117$ mod $243$.

Our method instead sees we are working mod $243$ and there are $2$ borrows so we treat $243$ as $3^{3+2}$.  This means our $n$ is $3$ and our $m$ is $2$. So we have $$\binom{2120211}{1202111}\equiv \binom{(21)(20)2}{(12)(02)1}\frac{\mathlarger{\binom{(20)21}{(02)11}}}{\mathlarger{\binom{(20)2}{(02)1}}}\frac{\mathlarger{\binom{211}{111}}}{\mathlarger{\binom{21}{11}}}\frac{\mathlarger{\binom{112}{110}}}{\mathlarger{\binom{11}{11}}}.$$  This converts to $90\frac{39}{3}\frac{26}{8}\frac{91}{1}$.  This is also congruent to $117$ mod $243$.

\bibliographystyle{plain}
\bibliography{bibliography.bib}

\begin{thebibliography}{1}

\bibitem{DW}
Kenneth~S. Davis and William~A. Webb.
\newblock Lucas' theorem for prime powers.
\newblock {\em European Journal of Combinatorics}, 11(3):229--233, 1990.

\bibitem{Kummer}
E.~E. Kummer.
\newblock Über die {Ergänzungssätze} zu den allgemeinen {Reciprocitätsgesetzen.}
\newblock {\em Journal für die reine und angewandte Mathematik}, 1852(44):93--146, 1852.

\end{thebibliography}

\end{document}